\documentclass[graybox]{svmult}

\usepackage{type1cm}        
\usepackage{makeidx}         
\usepackage{graphicx}        
\usepackage{multicol}        
\usepackage[bottom]{footmisc}

\usepackage{tikz}
\usepackage{array}

\usepackage{newtxtext}       %
\usepackage[varvw]{newtxmath}       

\makeindex             

\begin{document}

\title*{Multipreconditioning With Directional Sweeping Methods For High-Frequency Helmholtz Problems}
\titlerunning{Multipreconditioning With Directional Sweeping Methods For Helmholtz Problems}
\author{Niall Bootland\orcidID{0000-0002-3207-5395} and\\ Tyrone Rees\orcidID{0000-0003-0476-2259}}
\institute{Niall Bootland \at STFC Rutherford Appleton Laboratory, Harwell, UK, \email{niall.bootland@stfc.ac.uk}
\and Tyrone Rees \at STFC Rutherford Appleton Laboratory, Harwell, UK, \email{tyrone.rees@stfc.ac.uk}}
%
%
\maketitle

\abstract{We consider the use of multipreconditioning, which allows for multiple preconditioners to be applied in parallel, on high-frequency Helmholtz problems. Typical applications present challenging sparse linear systems which are complex non-Hermitian and, due to the pollution effect, either very large or else still large but under-resolved in terms of the physics. These factors make finding general purpose, efficient and scalable solvers difficult and no one approach has become the clear method of choice. In this work we take inspiration from domain decomposition strategies known as sweeping methods, which have gained notable interest for their ability to yield nearly-linear asymptotic complexity and which can also be favourable for high-frequency problems. While successful approaches exist, such as those based on higher-order interface conditions, perfectly matched layers (PMLs), or complex tracking of wave fronts, they can often be quite involved or tedious to implement. We investigate here the use of simple sweeping techniques applied in different directions which can then be incorporated in parallel into a multipreconditioned GMRES strategy. Preliminary numerical results on a two-dimensional benchmark problem will demonstrate the potential of this approach.}


\section{Introduction}
\label{sec:intro}

Developing effective and efficient solvers for the Helmholtz equation is a challenging task and presents an active area of ongoing research, especially for high-frequency problems where the wave number $k$ is large and so, due to the pollution effect \cite{BabuskaAndSauter-Pollution}, the linear systems to be solved typically become very large too. In this regime solutions are typically highly oscillatory and thus the underlying approximation should be able to capture this behaviour. Nonetheless, in real large-scale engineering applications the meshes used may under-resolve the problem to retain speed or limit the size of the system in order to feasibly solve it and this can exacerbate issues for the underlying solvers; see, for instance, \cite{BootlandEtAl-Comparison}.

Due to the size of the underlying linear systems, we must solve them using an iterative method, which then requires the use of a good preconditioner. 
Domain decomposition methods are natural candidates for this, as they combine the benefits of both direct and iterative solvers. 
In this work we consider an approach based sweeping methods, a promising family of techniques for high-frequency Helmholtz problems, and investigate the use of a multipreconditioning strategy which uses simple sweeps in different directions and allows for parallelisation in the preconditioner application (i.e., sweeping directions).

For a source term $f$, we wish to find the solution $u$ to the model free-space Helmholtz problem
\begin{subequations}
	\label{Helmholtz}
	\begin{align}
		\label{HelmholtzEquation}
		\Delta u + k^{2} u & = f & & \text{in } \Omega,\\
		\label{HelmholtzRobinBC}
		\frac{\partial u}{\partial n} - i k u & = 0 & & \text{on } \Gamma = \partial\Omega,
	\end{align}
\end{subequations}
where the wave number $k > 0$ is the quotient of the angular frequency $\omega$ and the wave speed $c$, namely $k = \omega / c$. 
We restrict ourselves to the free-space problem with impedance (or Robin) boundary conditions \eqref{HelmholtzRobinBC} for simplicity, 
although the solver methodology we consider can equally handle other boundary conditions. 
Indeed, boundary conditions which approximate a transparent condition, such as that in \eqref{HelmholtzRobinBC}, 
so as not to reflect waves back are fundamental in the approach of sweeping methods; see Section~\ref{sec:sweeping}. 
Afterwards we describe multipreconditioning in Section~\ref{sec:multiprec} and how we can used different sweeps to provide multiple preconditioners. 
Finally, some initial numerical results are presented in Section~\ref{sec:results}.

\section{Sweeping Methods For Helmholtz Problems}
\label{sec:sweeping}

Sweeping methods solve the Helmholtz via a sequence of smaller problems, 
namely a sequential decomposition into subdomains 
(the discrete equivalent being blocks in the matrix).
These methods sweep from the first subproblem to the last, 
transferring data along the way, and then sweep back to the first subproblem. 
At the discrete level this equates to using a block $LU$ factorisation, 
with the forward/backward sweeps being given by forward/backward substitution. 
In the continuous setting this requires the passing of boundary data 
from one subdomain to the next. 
Ideally this transmission condition (also called an interface condition) 
should let waves pass through the subdomains with no distortion or reflection 
from the artificial subdomain boundary 
and such a transparent condition is given by the so-called Dirichlet-to-Neumann (DtN) operator. 
The exact double sweep in this case then provides a nilpotent operator and yields the exact solution.

In practice the DtN operator is non-local and too expensive to work with and hence we must use approximate transparent boundary conditions (or equivalently approximate Schur complements in the discrete case within the block factorisation). Many techniques have been introduced to emulate transparent boundary conditions such as absorbing boundary conditions (ABCs), perfectly matched layers (PMLs) or the use of Pad\'e or rational interpolant approximation; see the survey paper \cite{GanderAndZhang}. The approach then becomes a preconditioner for an iterative method, such as Richardson iteration or a Krylov method.

Since a sweep is inherently sequential, much effort has been made to develop sweeping methods which are as efficient as possible in transmission conditions and method(s) of sweeping so as to reduce the overall computational cost. Nonetheless, this can lead to approaches which are quite complex and fiddly to implement or require in-depth modification of existing codes. Some approaches have already been developed to gain parallelism within sweeping methods, such as performing shorter sections of a sweep in parallel and decoupling the forward and backward sweeps \cite{VionAndGeuzaine-Parallel}, as well as separating the source term by subdomain and treating each separately to give an additive approach \cite{LiuAndYing-Additive,LengAndJu-Additive}. Here we will investigate the utility of incorporating parallelism through a multipreconditioning approach. In particular, we consider utilising double sweeps, for which a unified framework can be found in \cite{BouzianiEtAl-Framework}.

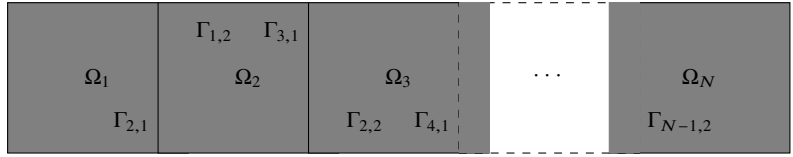
\begin{figure}[t]
	\begin{center}
		\begin{tikzpicture}[scale=2]
			\foreach \p in {0,1,2,4}
			{
				\fill[gray,opacity=0.4] (\p,0) -- (\p,1) -- (\p+1+0.2,1) -- (\p+1+0.2,0) -- cycle;
				\draw (\p,0) -- (\p,1) -- (\p+1+0.2,1) -- (\p+1+0.2,0) -- cycle;
			}
			\draw[dashed] (3,0) -- (3,1) -- (3+1+0.2,1) -- (3+1+0.2,0) -- cycle;
			\fill[gray,opacity=0.4] (3,0) -- (3,1) -- (3+0.2,1) -- (3+0.2,0) -- cycle;
			\fill[gray,opacity=0.4] (4,0) -- (4,1) -- (3+1+0.2,1) -- (3+1+0.2,0) -- cycle;
			\foreach \p in {1,2,3}
			{
				\draw (\p-0.4,0.5) node {$\Omega_{\p}$};
			}
			\draw (4+0.6,0.5) node {$\Omega_{N}$};
			\draw (3+0.6,0.5) node {$\cdots$};
			\draw (1,0.2) node[left] {$\Gamma_{2,1}$};
			\draw (1.2,0.8) node[right] {$\Gamma_{1,2}$};
			\draw (2.,0.8) node[left] {$\Gamma_{3,1}$};
			\draw (2.2,0.2) node[right] {$\Gamma_{2,2}$};
			\draw (3,0.2) node[left] {$\Gamma_{4,1}$};
			\draw (4.2,0.2) node[right] {$\Gamma_{N-1,2}$};
		\end{tikzpicture}
	\end{center}
	\caption{Sequential decomposition into $N$ subdomains $\Omega_{s}$, $1 \le s \le N$, with the internal subdomain boundaries $\Gamma_{s,1}$ (for $2 \le s \le N$) and $\Gamma_{s,2}$ (for $1 \le s \le N-1$) ordered sequentially as with the decomposition, here left-to-right in this illustration; the overlap regions are shaded darker.}
	\label{fig:OverlappingSeqDecomp}
\end{figure}

We stick to a simple sweep over $N$ sequential subdomains $\Omega_{s}$ for $1 \le s \le N$ with overlap. An illustration of the sequential decomposition is given in Figure~\ref{fig:OverlappingSeqDecomp}. For a double sweep with transmission conditions given by the operators $\mathcal{B}_{s,1}$ and $\mathcal{B}_{s,2}$ (approximating transparent conditions) on $\Gamma_{s,1}$ and $\Gamma_{s,2}$ respectively (see Figure~\ref{fig:OverlappingSeqDecomp}), we first solve the \emph{forward sweep} for $v_{s}$ in each subdomain given by
\begin{subequations}
	\label{HelmholtzForwardSweep}
	\begin{align}
		\label{HelmholtzEquationForward}
		\Delta v_{s} + k^{2} v_{s} & = f & & \text{in } \Omega_{i}, & & 1 \le s \le N,\\
		\label{HelmholtzBC1Forward}
		\mathcal{B}_{s,1}(v_{s}) & = \mathcal{B}_{s,1}(v_{s-1}) & & \text{on } \Gamma_{s,1}, & & 2 \le s \le N,\\
		\label{HelmholtzBC2Forward}
		\mathcal{B}_{s,2}(v_{s}) & = 0 & & \text{on } \Gamma_{s,2}, & & 1 \le s \le N-1,
	\end{align}
\end{subequations}
followed by the \emph{backward sweep} for $u_{s}$
\begin{subequations}
	\label{HelmholtzBackwardSweep}
	\begin{align}
		\label{HelmholtzEquationBackward}
		\Delta u_{s} + k^{2} u_{s} & = f & & \text{in } \Omega_{s}, & & 1 \le s \le N\\
		\label{HelmholtzBC1Backward}
		\mathcal{B}_{s,1}(u_{s}) & = \mathcal{B}_{s,1}(v_{s-1}) & & \text{on } \Gamma_{s,1}, & & 2 \le s \le N,\\
		\label{HelmholtzBC2Backward}
		\mathcal{B}_{s,2}(u_{s}) & = \mathcal{B}_{s,2}(u_{s+1}) & & \text{on } \Gamma_{s,2}, & & 1 \le s \le N-1,
	\end{align}
\end{subequations}
and where the problem boundary conditions, here \eqref{HelmholtzRobinBC}, are always maintained on $\Gamma = \partial\Omega$. We make use of a simple zeroth-order approximation to the transparent condition, as in \eqref{HelmholtzRobinBC}, and hence choose
\begin{align}
	\label{TransmissionCondition}
	\mathcal{B}_{s,j}(u) &= \frac{\partial u}{\partial n_{s,j}} - i k u, & & j \in \left\{1,2\right\},
\end{align}
where $n_{s,1}$ and $n_{s,2}$ are the outward normals from $\Omega_{s}$ along $\Gamma_{s,1}$ and $\Gamma_{s,2}$, respectively. Note that $u_{N} = v_{N}$ and so we only need one solve in the final subdomain $\Omega_{N}$. In the overlapping regions a partition of unity can be used to define a unique solution; in practice we will use the values from the most recent subdomain solve here.

We will investigate using these simple sweeps with decompositions in different directions which are then combined using multipreconditioning, as we now discuss. We note that multidirectional sweeping methods were also considered in \cite{DaiEtAl-Multidirectional} but for a non-overlapping checkerboard decomposition relying on high-order transmission conditions and, moreover, without the parallelism of multipreconditioning; instead, the preconditioner changes at each iteration of a flexible GMRES method.

\section{Multipreconditioned GMRES}
\label{sec:multiprec}

When using the standard (right-preconditioned) GMRES algorithm to solve $Ax=b$, with a preconditioner $M$, at iteration $k$ we find an $x_k$ within the so-called Krylov subspace $x_0 + \mathrm{span}\left\lbrace M^{-1}r_0, M^{-1}AM^{-1}r_0, \dots, (M^{-1}A)^{k-1}M^{-1}r_0\right\rbrace$ that minimises the Euclidean norm of the residual, where $r_0$ is the initial residual. Flexible GMRES allows us to change the preconditioner at each iteration, but still augments the search space with only one additional vector at each iteration, meaning it is inherently sequential.

The multipreconditioned GMRES method (MPGMRES) \cite{GreifEtAl-MultiplePreconditioners} allows the use of multiple preconditioners at the same time. If we have $t$ candidate preconditioners, say $M_1, \dots, M_t$, then \emph{complete MPGMRES} finds, at iteration $k$, the vector
\begin{align*}
	x_k = x_0 + \sum_{i=1}^{t} p_k^i(M_1^{-1}A, \dots, M_t^{-1}A)M_{i}^{-1}r_0
\end{align*}
that minimises the Euclidean norm of the residual. Here the $p_k^i(X_1, \dots, X_t)$ are multivariate $(k-1)$-degree polynomials in $t$ non-commuting variables (see \cite{GreifEtAl-MultiplePreconditioners} for further details). It is clear that this reduces to standard GMRES when used with only one preconditioner. The benefit of MPGMRES, however, is that it allows for cross terms between the preconditioners. This may be particularly advantageous when the preconditioners incorporate different physics or, as in our case, complementary directionality.

While complete MPGMRES is attractive, as it uses the complete space that can be constructed by applying multiple preconditioners simultaneously, it is, in general, infeasible, due to the exponential explosion in the dimension of the search space. A tractable algorithm is \emph{selective MPGMRES}, in which we augment the search space by $t$ additional vectors at each iteration. There are various ways this can be achieved, for example, if at step $k$ the search space was augmented by the columns of $V_k$, then we could apply all preconditioners to $V_k w_k$ for an appropriate weight vector $w_k$, or we could apply each preconditioner to a different column of $V_k$. Unless there is specific structure to exploit, all choices are heuristic, and software implementing this technique will offer a range of selective MPGMRES schemes.

We remark that at the $k$th iteration of selective MPGMRES, applied using $t$ preconditioners, we require $t$ matrix--vector products and $t$ preconditioner solves along with $(k-\frac{1}{2})t^2+\frac{3}{2}t$ inner products. In contrast, $t$ iterations of a flexible GMRES method cycling through $t$ preconditioners would cost the same number of matrix--vector products and preconditioner solves, but now $\frac{k}{2}t^2 + (\frac{k}{2}+1)t$ inner products. However, the preconditioner solve is typically the most expensive part of a Krylov subspace method and we can perform this part fully in parallel within MPGMRES while, at the same time, typically converging faster due to the enriched search space.

\section{Numerical Results}
\label{sec:results}

We consider a simple benchmark problem of a square domain $\Omega = (0,1)^2$ with a point source in the centre given by
\begin{align*}
	f(x,y) = 3\times10^{4}\exp\left(-200k\left((x-0.5)^2 + (y-0.5)^2\right)\right).
\end{align*}
To discretise the problem we make use of the standard five-point finite different approximation on a regular Cartesian grid. In forming the preconditioners we split the square into $N$ regular strip-wise subdomains, either in $x$ or $y$, and add one layer of mesh points for every internal boundary to give an overlapping decomposition; the overlap between neighbouring subdomains is then of width $2h$ where $h$ is the mesh spacing. One preconditioner uses a decomposition in $x$ (left-to-right) and one uses a decomposition in $y$ (bottom-to-top). We use \texttt{MATLAB} for implementation, using \texttt{$\backslash$} for subdomain problems, along with the MPGMRES code of the second author\footnote{https://uk.mathworks.com/matlabcentral/fileexchange/34562-multi-preconditioned-gmres}. A relative residual tolerance of $10^{-6}$ is used throughout for the MPGMRES solve.

We first compare different combinations of sweeping strategies. We denote by LRL a double sweep in the horizontal direction, namely from the left to the right and back to the left again. Similarly, BTB denotes a double sweep in the vertical direction and we use LR, RL, BT and TB to refer to single (forward) sweeps going right, left, up and down, respectively. When combining in the multipreconditioning step we use a $+$ to denote this, for instance LRL+BTB refers to the ordered combining of LRL and BTB double sweeps. Where only two preconditioners are combined we apply the preconditioner only to the previous vector corresponding to the other preconditioner. For more than two preconditioners we apply the preconditioner to the sum of all vectors at the previous iteration.

In Table~\ref{tab:comp_precs} we consider four sweeping combinations, a single LRL double sweep (no multipreconditioning), the combination of double sweeps LRL+BTB and two similar approaches only using forward sweeps, namely LR+RL and LR+BT+RL+TB. Note that in the latter two approaches each preconditioner is roughly half as expensive as the former (unless the double sweeps are parallelised by decoupling forward and backward sweeps as in \cite{VionAndGeuzaine-Parallel}). We observe similar results for each $k$, the most notable of which is that LRL+BTB takes half the iterations compared to LRL. This suggests that multipreconditioning can be highly beneficial and a parallel implementation would reduce the computation time by a factor of two in this case. Splitting the forward and backward sweeps to be two separate forward sweeps increases the iteration counts, as expected, but noticeably not by double in the case of LR+RL compared to LRL where twice as much work is required in the preconditioner application. However, it is not as successful when combining different coordinate directions and we see no benefit in LR+BT+RL+TB since here the iteration counts double, nullifying any benefit from the cheaper preconditioner applications and, moreover, increasing the number of inner products required within MPGMRES. As such, the most promising combination is in double sweeps along different directions, as given by LRL+BTB and we now focus on this approach.

\begin{table}[!t]
	\caption{Iteration counts with different combinations of sweeping preconditioners for $N=8$. The wavenumber $k$ is varied for a fixed mesh size of $h = 2^{-9}$.}
	\label{tab:comp_precs}
	\centering
	\begin{tabular}{p{1cm}|>{\centering}p{2.4cm}>{\centering}p{2.4cm}>{\centering}p{2.4cm}>{\centering\arraybackslash}p{2.4cm}}
		\hline\noalign{\smallskip}
		$k$ & LRL & LRL+BTB & LR+RL & LR+BT+RL+TB \\
		\noalign{\smallskip}\svhline\noalign{\smallskip}
		50 & 20 & 9 & 31 & 18 \\
		100 & 18 & 8 & 29 & 18 \\
		200 & 17 & 9 & 31 & 18 \\
		\noalign{\smallskip}\hline\noalign{\smallskip}
	\end{tabular}
\end{table}

We now consider more carefully the performance of LRL and LRL+BTB as the wavenumber $k$ increases. In Table~\ref{tab:vary_k} we see that there is a regime where iteration counts remain nearly constant as $k$ increases, for the fixed $h = 2^{-9}$ mesh used this extends to around $k=320$. However, beyond this the number of points per wavelength (ppwl) decreases below 10 and we start to see some degradation in iteration counts as $k$ continues to increase. This degradation is slightly more mild for LRL compared to LRL+BTB but the latter still yields the lowest iteration counts up to our maximum wavenumber $k=640$, which is approximately 5 ppwl. Note that while iteration counts are not fully robust, we are still able to solve the problem in relatively few iterations even in the most challenging cases (cf. \cite{BootlandEtAl-Comparison} where some methods are failing at 5 ppwl, albeit on larger and more realistic test cases).

\begin{table}[!t]
	\caption{Iteration counts with double sweep approaches for $N=8$ and varying wavenumber $k$. A fixed mesh size of $h = 2^{-9}$ is used throughout.}
	\label{tab:vary_k}
	\centering
	\begin{tabular}{l|>{\centering}p{0.5cm}>{\centering}p{0.5cm}>{\centering}p{0.5cm}>{\centering}p{0.5cm}>{\centering}p{0.5cm}>{\centering}p{0.5cm}>{\centering}p{0.5cm}>{\centering}p{0.5cm}|>{\centering}p{0.5cm}>{\centering}p{0.5cm}>{\centering}p{0.5cm}>{\centering}p{0.5cm}>{\centering}p{0.5cm}>{\centering}p{0.5cm}>{\centering}p{0.5cm}>{\centering\arraybackslash}p{0.5cm}}
		\hline\noalign{\smallskip}
		$k$ & 40 & 80 & 120 & 160 & 200 & 240 & 280 & 320 & 360 & 400 & 440 & 480 & 520 & 560 & 600 & 640 \\
		\noalign{\smallskip}\svhline\noalign{\smallskip}
		LRL & 20 & 18 & 18 & 17 & 17 & 18 & 18 & 18 & 19 & 20 & 21 & 23 & 23 & 26 & 27 & 29 \\
		LRL+BTB & 9 & 8 & 8 & 9 & 9 & 10 & 10 & 11 & 13 & 13 & 14 & 16 & 17 & 19 & 22 & 27 \\
		\noalign{\smallskip}\hline\noalign{\smallskip}
	\end{tabular}
\end{table}

We also consider decoupling concerns of ppwl by varying the mesh with $k$ so as to ensure a discretisation using 10 ppwl. Results in Table~\ref{tab:vary_k_10ppwl} now show that the dependence on $k$ is worse for LRL while there is now only a very mild increase in iteration counts for LRL+BTB, which is able to solve the problem with $k=643.4$ in just 12 iterations. Overall this presents very positive results for a multipreconditioned sweeping approach in the high frequency (large $k$) regime, especially when the mesh is sufficiently refined (e.g., using 10 ppwl).

\begin{table}[!t]
	\caption{Iteration counts with double sweep approaches for $N=8$ and varying wavenumber $k$ with the mesh spacing $h$ chosen to ensure a fixed 10 ppwl.}
	\label{tab:vary_k_10ppwl}
	\centering
	\begin{tabular}{l|>{\centering}p{0.8cm}>{\centering}p{0.8cm}>{\centering}p{0.8cm}>{\centering}p{0.8cm}>{\centering}p{0.8cm}>{\centering\arraybackslash}p{0.8cm}}
		\hline\noalign{\smallskip}
		$h$ & $2^{-5}$ & $2^{-6}$ & $2^{-7}$ & $2^{-8}$ & $2^{-9}$ & $2^{-10}$ \\
		$k$ & 20.1 & 40.2 & 80.4 & 160.8 & 321.7 & 643.4 \\
		\noalign{\smallskip}\svhline\noalign{\smallskip}
		LRL & 8 & 10 & 12 & 15 & 19 & 22 \\
		LRL+BTB & 8 & 10 & 10 & 11 & 11 & 12 \\
		\noalign{\smallskip}\hline\noalign{\smallskip}
	\end{tabular}
\end{table}

To investigate the effect of the mesh resolution for a fixed wavenumber, in Table~\ref{tab:vary_h} we vary $h$ for a fixed $k=50$. We observe almost constant iteration counts for LRL+BTB for $h=2^{-6}$ to $h=2^{-10}$ (with over a million unknowns) while iteration counts more than double over the same range when just using LRL. Note that for $h=2^{-5}$ we have just 4 ppwl and the discretisation is unable to capture the true solution and so both approaches require more iterations to converge; this represents a typical case with iterative methods where the solver will struggle if the discrete problem is not a meaningful representation of the underlying continuous problem.

\begin{table}[!t]
	\caption{Iteration counts with double sweep approaches for $N=8$ and varying the mesh spacing $h$ for a fixed wavenumber $k = 50$.}
	\label{tab:vary_h}
	\centering
	\begin{tabular}{l|>{\centering}p{0.8cm}>{\centering}p{0.8cm}>{\centering}p{0.8cm}>{\centering}p{0.8cm}>{\centering}p{0.8cm}>{\centering\arraybackslash}p{0.8cm}}
		\hline\noalign{\smallskip}
		$h$ & $2^{-5}$ & $2^{-6}$ & $2^{-7}$ & $2^{-8}$ & $2^{-9}$ & $2^{-10}$ \\
		\noalign{\smallskip}\svhline\noalign{\smallskip}
		LRL & 21 & 10 & 12 & 15 & 20 & 24 \\
		LRL+BTB & 25 & 10 & 8 & 8 & 9 & 10 \\
		\noalign{\smallskip}\hline\noalign{\smallskip}
	\end{tabular}
\end{table}

Finally, in Table~\ref{tab:vary_N_time} we vary 
the number of subdomains used for each sweeping method. 
We see that, as expected, iteration counts increase with $N$, 
nonetheless this increase is more noticeable for LRL compared with LRL+BTB. 
For $k=50$, in going from $N=4$ to $N=64$ we note that the iteration counts for 
LRL more than triple while for LRL+BTB they only double. 
Further, by $N=64$ the multipreconditioned approach requires just a quarter 
of the iterations needed by the LRL double sweep. 
The trend when $k=100$ is slightly poorer for LRL+BRB but the approach is still highly favourable, 
now requiring just over a third of the iterations of LRL for $N=64$.

In Table~\ref{tab:vary_N_time} we also include the wall-clock time to run these tests 
in a serial implementation using Matlab R2024a on a laptop with an
Intel i7-1270P processor running at 2.20 GHz, with 32GB of RAM. We highlight that, even 
with a sequential implementation, using sweeping with MPGMRES can give a significant speedup 
when computing with more subdomains.

\begin{table}[!t]
	\caption{Iteration counts (and wall-clock time in seconds) with double sweep approaches for varying $N$ with fixed mesh spacing $h = 2^{-9}$ and wavenumbers $k = 50$ and $k=100$.
	Note that the times are for a serial implementation, where the preconditioners are not applied in parallel.}
	\label{tab:vary_N_time}
	\centering
	\begin{tabular}{l|>{\centering}p{1.8cm}>{\centering}p{1.8cm}>{\centering}p{1.8cm}>{\centering}p{1.8cm}>{\centering\arraybackslash}p{1.8cm}}
	\hline\noalign{\smallskip}
	& \multicolumn{5}{c}{$k=50$} \\
	$N$ & $4$ & $8$ & $16$ & $32$ & $64$\\
	\noalign{\smallskip}\svhline\noalign{\smallskip}
	LRL & 17 (8.3s) & 20 (8.6s) & 24 (9.7s) & 37 (15s) & 60 (25s)\\ 
	LRL+BTB & 8 (8.9s) & 9 (8.8s) & 10 (8.8s) & 13 (13s) & 15 (14s)\\ 
	\hline\noalign{\smallskip}
	& \multicolumn{5}{c}{$k=100$} \\
	$N$ & $4$ & $8$ & $16$ & $32$ & $64$\\
	\noalign{\smallskip}\svhline\noalign{\smallskip}
	LRL & 16 (8.5s)& 18 (11s)& 22 (9.3s)& 33 (13s)& 52 (24s)\\ 
	 LRL+BTB & 7 (9.1s)& 8 (8.4s)& 10 (10s)& 13 (11s)& 19 (16s)\\ 
	\noalign{\smallskip}\svhline\noalign{\smallskip}
	\end{tabular}
\end{table}

The results presented here are highly promising 
and suggest that incorporating sweeps in different directions through multipreconditioning 
can be very effective to reduce iteration counts. 
We observe wall-clock times competitive with a single preconditioner when run in serial, 
and envisage that running the preconditioners in parallel using MPI or OpenMP would 
show a more significant speed-up in the computation time. 
Multipreconditioning provides a novel way to help parallelise sweeping methods, 
which by themselves have an inherently sequential nature. 
There are many avenues to further explore this approach, 
as well as in combining it with other techniques to parallelise the overall solver. 
More realistic test problems must be considered as well as larger scale tests. 
Further, a more efficient and parallel implementation is required 
to reveal the wall clock savings such an approach can bring; while the application of
the preconditioners is trivially parallelizable, the effect of the inherent overheads in
a parallel environment plus communication between the processors, remains to be seen. 
We note a parallel version of MPGMRES is available as \texttt{HSL\_MI29} in the HSL library\footnote{https://www.hsl.rl.ac.uk/catalogue/index.html}.

In terms of the underlying sweeping methods, we have considered one of the simplest transmission conditions in \eqref{TransmissionCondition}, allowing easy implementation. It would be pertinent to also investigate more accurate transmission conditions, such as through higher order ABCs or PMLs, and verify the effectiveness of multipreconditioning in these cases. In the present scheme, since there is no second level, the method does not scale with $N$ and so relatively few larger subdomains is best but the subproblems can then equally be solved by the same approach and further decomposed recursively, say in a perpendicular direction, to give more parallelism and reduce the burden on direct solvers. Preliminary results (not shown here) suggest the subproblems can be solved to low accuracy without much deterioration in the overall solver and this should be explored further. Moreover, in 3D there are three coordinate directions for sweeping which may further enhance parallelism and the effectiveness of our approach for such problems and it will be interesting to study this. Further, diagonal sweeps on a Cartesian decomposition can also be considered, as in \cite{DaiEtAl-Multidirectional}. Finally, we note that our multipreconditioned sweeping approach could also be applied to the Maxwell equations, where higher order ABCs are trickier to deal with.

\bibliographystyle{spmpsci}
\bibliography{refs}

\end{document}